\documentclass[12pt]{article}
\usepackage[latin1]{inputenc}
\usepackage{amsmath}
\usepackage{amsfonts}
\usepackage{amssymb}
\usepackage{amsthm}

\usepackage{mathrsfs} 

\usepackage[alphabetic,lite]{amsrefs} 

\usepackage{fullpage}
\usepackage{setspace}
\usepackage{hyperref}
\usepackage{color}
\usepackage{enumerate} 

\usepackage[all,cmtip]{xy} 

\usepackage{fancyhdr}
\newtheorem{Thm}{Theorem}[section]
\newtheorem{Lemma}[Thm]{Lemma}
\newtheorem{Cor}[Thm]{Corollary}

\newtheorem{Prop}[Thm]{Proposition}

\theoremstyle{definition}

\newtheorem{Rem}[Thm]{Remark}

\newcommand{\<}{\left\langle}
\renewcommand{\>}{\right\rangle}

\newcommand{\gl}{\lambda}

\newcommand{\al}{\alpha}

\newcommand{\lam}{\lambda}

\newcommand{\bC}{\mathbb{C}}

\newcommand{\bF}{\mathbb{F}}
\newcommand{\bG}{\mathbb{G}}

\newcommand{\bQ}{\mathbb{Q}}
\newcommand{\bR}{\mathbb{R}}

\newcommand{\bZ}{\mathbb{Z}}

\newcommand{\cG}{\mathcal{G}}

\newcommand{\cO}{\mathcal{O}}

\newcommand{\ff}{\mathfrak{f}}
\newcommand{\fg}{\mathfrak{g}}

\newcommand{\fs}{\mathfrak{s}}

\newcommand{\fF}{\ov \bF_p}

\newcommand{\sA}{\mathscr{A}}
\newcommand{\sB}{\mathscr{B}}

\newcommand{\sG}{\mathscr{G}}
\newcommand{\sH}{\mathscr{H}}

\newcommand{\sT}{\mathscr{T}}

\newcommand{\sV}{\mathscr{V}}

\DeclareMathAlphabet{\mathpzc}{OT1}{pzc}{m}{it}

\newcommand{\ra}{\rightarrow}

\newcommand{\ov}{\overline}

\DeclareMathOperator{\Spec}{Spec}		
\DeclareMathOperator{\Hom}{Hom}			
\DeclareMathOperator{\GL}{GL}		  	
\DeclareMathOperator{\SL}{SL}		  	
\DeclareMathOperator{\SO}{SO} 			
\DeclareMathOperator{\Gal}{Gal}			
\DeclareMathOperator{\Aut}{Aut}			
\DeclareMathOperator{\Lie}{Lie}  	  
\DeclareMathOperator{\ind}{ind}  	  
\DeclareMathOperator{\Sym}{Sym}  	  
\DeclareMathOperator{\val}{\mathsf{v}}  	  
\newcommand{\Waff}{W_{\mathrm{aff}}}  	  
\DeclareMathOperator{\disc}{disc}  
\newcommand{\ba}{\begin{aligned}}
	\newcommand{\ea}{\end{aligned}}
\newcommand{\beqn}{\begin{eqnarray}}
	\newcommand{\eeqn}{\end{eqnarray}}
\newcommand{\beqns}{\begin{eqnarray*}}
	\newcommand{\eeqns}{\end{eqnarray*}}
\newcommand{\benum}{\begin{enumerate}}
	\newcommand{\eenum}{\end{enumerate}}

\newcommand{\ZG}{\sH}       
\newcommand{\ZV}{\sV}

\newcommand{\ZZG}{\sG}
\newcommand{\ZZT}{\sT}

\newcommand{\RP}{\mathbf{G}} 
\newcommand{\PV}{\mathbf{V}} 
\newcommand{\RV}{\PV}

\newcommand{\mmu}{\boldsymbol\mu}
\DeclareMathOperator{\Conjug}{Inn} 
\newcommand{\Kal}{\ov \bQ_p} 
\newcommand{\Oal}{\mathbf O} 
\newcommand{\Root}{\Psi} 
\newcommand{\Gk}{\cG} 
\newcommand{\Weil}{\mathcal W} 
\DeclareMathOperator{\Stab}{Stab} 

\usepackage[margin=1in]{geometry}

\AtEndDocument{\bigskip{\small%
		Jessica Fintzen ~~ \texttt{fintzen@math.harvard.edu} \\		
		 Department of Mathematics, Harvard University, One Oxford Street, Cambridge, MA 02138, USA   \\		 
		\textit{Current address}: Department of Mathematics, University of Michigan, 2074 East Hall, 530 Church St, Ann Arbor, MI 48109, USA \par
		\addvspace{\medskipamount}
	Beth Romano ~~\texttt{romanob@bc.edu}\\ 
	Department of Mathematics, Boston College, 140 Commonwealth Avenue, Chestnut Hill, MA 02467, USA\\ 
		\textit{Current address}: Department of Pure Mathematics and Mathematical Statistics, University of Cambridge, Wilberforce Road, Cambridge, CB3 0WB, UK
		\par
	}}

\begin{document}
	\author{Jessica Fintzen and Beth Romano}
	\title{Stable vectors in Moy--Prasad filtrations}
	\date{}
	\maketitle

	\begin{abstract} 
	Let $k$ be a finite extension of $\bQ_p$, let ${\Gk}$ be an absolutely simple split reductive group over $k$, and let $K$ be a maximal unramified extension of $k$. 
	To each point in the Bruhat--Tits building of $\Gk_K$, Moy and Prasad have attached a filtration of $\Gk(K)$ by bounded subgroups.
	In this paper we give necessary and sufficient conditions  for the dual of the first Moy--Prasad filtration quotient 
	to contain stable vectors for the action of the  reductive quotient.  
			
	Our work extends earlier results by Reeder and Yu, who gave a classification in the case when $p$ is sufficiently large. By passing to a finite unramified extension of $k$ if necessary, we obtain new supercuspidal representations of ${\Gk}(k)$. 
	\end{abstract}

{
\renewcommand{\thefootnote}{}  
	\footnotetext{The first authors's research was partially supported by the Studienstiftung des deutschen Volkes.}
}

	\tableofcontents

\section{Introduction}
	
Let $k$ be a finite extension of $\bQ_p$, {and let} ${\Gk}$ {be} an absolutely simple split reductive group over $k$. {In \cite{Yu},} Yu gave a construction of supercuspidal representations of ${\Gk}(k)$ with complex coefficients, generalizing earlier work of Adler (\cite{Adler}). {F}or large primes $p$, Yu's construction yields all possible supercuspidal representations {(\cite{Kim}).} However, the construction does not capture all supercuspidal representations {when $p$ is small}. 
In \cite{ReederYu}, Reeder and Yu gave a new construction of supercuspidal representations of lowest positive depth, which they call \textit{epipelagic} representations. {Their construction is} uniform for all $p$, but requires as input stable vectors in a certain representation coming from a Moy--Prasad filtration, which we will describe below.  For $p$ large enough, necessary and sufficient conditions for the existence of stable vectors were given in \cite{ReederYu}. 
For small $p$, the existence of stable vectors, hence also of epipelagic representations,  was previously unknown, except in the simple case of a torus action (\cite[Section~9]{GrossReeder}, \cite[Section~2.6]{ReederYu}). 

In this paper, we show that the criterion of \cite{ReederYu} for the existence of stable vectors is valid for all primes $p$, uniformly. Thus \cite{ReederYu} can be used to construct supercuspidal representations uniformly for all $p$:
\theoremstyle{plain}
\newtheorem*{Corintro}{Corollary \ref{Cor-main}}
\begin{Corintro} Let ${\Gk}$ be an absolutely simple split reductive group over {the local field} $k$. Then for each {$m$} satisfying the conditions of Theorem \ref{main-thm}, there exists a finite unramified extension $k'$ of $k$ such that one can implement the construction of \cite{ReederYu} to form supercuspidal epipelagic representations of ${\Gk}(k')$.
\end{Corintro}

For small $p$, we obtain previously unknown representations (see Remark \ref{Rem-main}).

Let $K$ be a maximal unramified extension of $k$ with ring of integers $\cO$ and residue field ${\fF}$. In order to describe the input for the construction of \cite{ReederYu}, we choose a rational point $x$ in the Bruhat--Tits building $\sB(G,K)$ of $G:=\Gk_K$. Moy and Prasad defined in \cite{MP1, MP2} a filtration of $G(K)$ {by subgroups:}
\begin{eqnarray*}
	G(K)_{x,0} \supsetneq G(K)_{x, r_1} \supsetneq G(K)_{x, r_2} \supsetneq \hdots  
\end{eqnarray*}
{where $r_1, r_2, ...$ are rational numbers depending on $x$.}
The quotient $G(K)_{x,0}/G(K)_{x,r_1}$ can be identified with the $\fF$-points of a reductive group $\RP_x${. The quotient $\RV_x :=$} $G(K)_{x,r_1}/G(K)_{x,r_2}$ {forms} an $\fF$-vector space {on which $\RP_x$ has an action induced by conjugation}. 
A vector in {the dual space} $\check \PV_x$ is said to be \textit{stable} {under the action of $\RP_x$} if its orbit is Zariski closed and its stabilizer in $\RP_x$ is finite. The existence of such a stable vector {over $k$} is the only requirement for Reeder--Yu's construction of epipelagic representations. 

Our main result{, which is independent of $p$,} is the {classification of all such $x$ for which the representation of $\RP_x$ on $\check \RV_x$ contains stable vectors}.
{To state the theorem, we fix an apartment} $\sA$ in $\sB(G,K)$ {and} a hyperspecial  point ${x_0 \in} \sA$. {We use} $\check \rho$ to denote the half-sum of a set of positive coroots.

\theoremstyle{plain} 
\newtheorem*{Thmintro}{Theorem \ref{main-thm}}
 
\begin{Thmintro} Let $x \in \sA$. Then the representation $\check \RV_{x}$ contains stable vectors under the action of $\RP_x$ if and only if $x$ is conjugate under the affine Weyl group $\Waff$ to $x_0 + \frac{1}{m}\check\rho$, where $m$ is the order of an elliptic, $\bZ$-regular element in the Weyl group $W$ of $G$.
\end{Thmintro}

For the definitions of elliptic and $\bZ$-regular see Section \ref{section-Vinberg}. We also obtain a similar result about the classification for semistable vectors (Proposition \ref{Prop-semistable}).

These results follow from a more general analysis of the relations between (semi)stable vectors in the special and geometric generic fibers of representations over {$\cO$}. More explicitly, {let} $\ZZG$ {be} a split reductive group over ${\cO}$ acting on a free ${\cO}$-module $\ZV$, {and let ${\Kal}$ be an algebraic closure of $K$}. We show that the representation of $\ZZG_{\Kal}$ on $\ZV_{\Kal}$ contains stable vectors if and only if the representation of $\ZZG_{\fF}$ on $\ZV_{\fF}$ contains stable vectors (Corollary \ref{Cor-stable}). The same statement is true with ``stable" replaced by ``semistable." Thus we can transfer results about the existence of (semi)stable vectors in characteristic zero to arbitrary positive characteristics. This is a key step in the proof of Theorem \ref{main-thm}.

\textbf{Structure of the paper.} In Section \ref{section-lifting-descending} we first recall basic definitions and properties related to stability. We then show that the special fiber of a representation of a split reductive group over ${\cO}$ has (semi)stable vectors if and only if the {geometric} generic fiber of this representations admits (semi)stable vectors. In Section \ref{section-Moy--Prasad}, we review Moy--Prasad filtrations and define a split reductive group over ${\cO}$ acting on a free ${\cO}$-module such that the special fiber of  this action corresponds to the action of $\RP_x$ on $\RV_x$. The generic fiber is isomorphic to a representation coming from Vinberg theory of graded Lie algebras, which we review in Section \ref{section-Vinberg}.
{We finish the proof of Theorem \ref{main-thm} u}sing results from {\cite{RLYG} on} Vinberg theory in characteristic zero. In the {final} section, Section \ref{section-G2}, we {give} an example{, classifying} all stable vectors in $\check \RV_x$ for the case of $G=G_2$ and $x=x_0+\frac{1}{2}\check \rho$. We then determine the Langlands parameter corresponding to the unique representation of $G_2(\bQ_2)$ obtained using this choice of $x$ in the construction of \cite{ReederYu}.

\textbf{Conventions and notation.} Throughout the paper, we use the following conventions. If $A$ is a ring, $H$ a group scheme over $\Spec A$, and $R$ an $A$-algebra, we write $H(R)$ for the $\Spec R$-points of $H$ and $H_R$ for $H \times_{\Spec A} \Spec R$. For a free ${A}$-module $V$, we write $V_R$ to denote the free $R$-module $V \otimes_{A} R$. If {$M$ is an $A$-module and} there is no danger of confusion, we may also denote by $M$ the scheme corresponding to the functor of points $R \mapsto M \otimes_A R$ for any $A$-algebra $R$.
	
	Throughout the paper, we consider reductive groups to be connected. 
	
	We maintain the following notation. Let $k$ be a finite field extension of $\bQ_p$, and let $K$ be a maximal unramified extension of $k$. We fix a discrete valuation $\val$ on $K$ with image $\bZ \cup \{\infty\}$. Let $\cO$ denote the ring of integers of $K$, and let $\varpi$ be a uniformizer of $\cO$. We fix an isomorphism between the residue field $\cO/\varpi \cO$ and ${\overline{\bF}_{{p}}}$.  We fix an algebraic closure ${\Kal}$ of $K$ with {ring of integers} $\Oal$.
	
	\textbf{Acknowledgements.} 
 The authors would like to thank their PhD advisors, Benedict Gross and Mark Reeder, for their support. 	Many of the results in this paper are adapted from the PhD dissertations of the authors. They also thank Mark Reeder for carefully reading an initial draft of this paper. The second author would like to thank Maksym Fedorchuk for a helpful conversation about the example in Section \ref{section-G2}. The authors also appreciate the referee's detailed review and feedback.

	\section{Lifting and descent of stable vectors} \label{section-lifting-descending} \label{Section-lifting}
	
	Before describing our results about {the} lifting and {descent of} (semi)stable vectors, we will recall basic definitions and properties about stability, and introduce some notation. 
	
	If $G$ is a reductive group over an algebraically closed field $F$, and $V$ is a finite-dimensional algebraic representation of $G$,
	a vector $v \in V$ is called \textit{semistable} if its orbit under $G(F)$ does not contain zero in its closure under the Zariski topology on $V$, or, equivalently, if there exists a non-constant $G(F)$-invariant homogeneous polynomial $f$ on $V$ such that $f(v) \neq 0$. A vector $v \in V$ is called \textit{stable} if its orbit is closed and its stabilizer in $G$ is finite.  In particlar, a stable vector is automatically semistable.

	In this paper we are interested in (semi)stable vectors in the following setting.
	Let $\ZZG$ be a split reductive group scheme over $\cO$ with split maximal torus $\ZZT$.
	Let $\ZV$ be a free finite-dimensional $\cO$-module, and let $\ZZG \to \GL(\ZV)$ be a representation of $\ZZG$, where, by abuse of notation, $\GL(\ZV)$ denotes the $\cO$-group scheme whose functor of points is given by $A \mapsto \GL(\ZV_A)$. By base change, we obtain representations $\ZZG_{\Kal} \to \GL(\ZV_{\Kal})$ and $\ZZG_{\fF} \to \GL(\ZV_{\fF})$. We will investigate the relationship between the geometric invariant theory of these two representations.
	To do this, we will make frequent use of the fact that since $\ZZG$ is smooth, $\ZZG(\cO)$ (and thus $\ZZG(\Oal)$) surjects onto $\ZZG({\fF})$.  
	We denote the image of an element $g \in \ZZG(\Oal)$ in $\ZZG({\fF})$ by $\overline g$.
	Similarly, we denote the surjection $\ZV_{\Oal} \twoheadrightarrow \ZV_{\fF}$ by $v \mapsto \overline v$. 

	Our other main tool is the Hilbert--Mumford Criterion (\cite[1.1]{Mumford}, which is based on \cite[Theorem~2.1]{Mumfordbook}). Recall that for any $\cO$-algebra  $A$ and one-parameter subgroup $\lam: \mathbb{G}_m \to \ZZG_A$, we obtain a weight decomposition $\ZV_A = \oplus_{i \in \bZ} (\ZV_A)_i$, where $\lam (t) \cdot v = t^i v$ for all $v \in (\ZV_A)_i$, $t\in\bG_m(A)\simeq A^\times$. Given a vector $v \in \ZV_A$, we let $I_A(\lam, v)$ denote the set of negative weights for $v$: if we write $v$ as a sum $v = \sum v_i$ such that $v_i \in (\ZV_A)_i$, then $I_A(\lam, v) = \{ i < 0 \mid v_i \neq 0 \}$. If $A$ is an algebraically closed field, e.g. ${\Kal}$ or ${\fF}$, the Hilbert--Mumford Criterion says that a vector $v \in \ZV_A$ is stable under $\ZZG_A$ if and only if for every nontrivial one-parameter subgroup $\lam: \mathbb{G}_m \to \ZZG_A$, {we have} $I_A(\lam, v) \neq \varnothing$. 
	We note for later use that if $g \in \ZZG(A)$, then $I_A(g\lam g^{-1}, v) = I_A(\lam, {g^{-1}} v)$.

	\begin{Lemma} \label{Lemma-1}
		If the representation $\ZV_{\Kal}$ contains stable vectors for the action of $\ZZG_{\Kal}$, then the representation $\ZV_{\fF}$ contains stable vectors for the action of $\ZZG_{\fF}$.
	\end{Lemma}
	
	\textbf{Proof.}
	Suppose the set $(\ZV_{\Kal})_s$ of stable vectors in $\ZV_{\Kal}$ is non-empty.
	Since $(\ZV_{{\Kal}})_s$ is open (see {\cite[1.4, p. 37]{Mumfordbook}}), there exists a nonzero polynomial $P$ on $\ZV_{\Kal}$ such that if $v \notin (\ZV_{\Kal})_s$, then $P(v) = 0$. Note that since $0 \in \ZV{_{\ov \bQ_p}}$ is not stable, we have $P(0) = 0$. Choosing a basis for $\ZV$ (and thus for $\ZV_{\Kal}$), we may identify $P$ with an element of ${\Kal} [x_1, ..., x_n]$.
	We can assume without loss of generality that the minimum of the valuations of the coefficients of $P$ is zero.  
	Let $\ov P$ be the image of $P$ under the reduction map $\Oal[x_1, \hdots, x_n] \ra {\fF}[x_1, \hdots, x_n] $. Then $\ov P \neq 0$ by choice of the coefficients, and $\overline P$ is not constant, because $P(0)=0$. Thus there exists $\ov v \in \ZV_{{\fF}}$ such that $\ov P(\ov v) \neq 0$. We claim that $\ov v$ is a stable vector under the action of $\ZZG_{\fF}$.
	Assume it is not. Then there exists a nontrivial one-parameter subgroup $\overline \lam: \mathbb{G}_m \to \ZZG_{\fF}$ such that $I_{\fF}(\ov \lam, \ov v) = \varnothing$. Choose an element ${\ov g} \in \ZZG({\fF})$ such that $\bar \lam_1 := \bar g\bar \lam \bar g^{-1}$ has image in $\ZZT({\fF})$. We have $I_{\fF}(\bar \lam_1, \bar g \bar v) = I_{\fF}(\bar \lam, {\ov v}) = \varnothing$.
	Using the identification of the cocharacters of $\ZZT_{\fF}$ and the cocharacters of $\ZZT_\Oal$, we obtain a lift $\lam_1: \mathbb{G}_m \to \ZZT_{\Kal} \hookrightarrow \ZZG_{\Kal}$ of $\bar \lam_1$. We also lift ${\ov v}$ to an element $v \in \ZV$ and lift ${\ov g}$ to an element $g \in \ZZG(\cO)$. Since $P(v) \neq 0$, $v$ is stable, and so $w := {g} v$ is also stable. Using the decomposition $\ZV=\bigoplus_{i \in \bZ} \ZV_i$ coming from $\lambda_1$, we write $w = w_+ + w_-$, where $w_+ \in \oplus_{m \in \bZ_{{\geq} 0}} \ZV_m $ and $w_- \in \oplus_{m \in \bZ_{<0}} \ZV_m$. Since $I_{\fF}(\bar \lam_1, {\ov w}) =I_{\fF}(\bar \lam_1, \bar g \bar v) = \varnothing$, we have $\overline{w}_- = 0$, and so $w_+$  is also a lift of ${\ov w}$ in $\ZV$. Note that $I_{\Kal}(\lam_1, w_+) = \varnothing$, so $w_+$ is not stable as an element of $\ZV_{\Kal}$. 
	Consider $ v' := g^{-1}(w_+)$. Since $\overline{g^{-1}(w_+)} = \overline{g}^{-1}{\ov w}= \overline{g}^{-1}\bar g \bar v$, we see that 
	$v'$ is a lift of $\bar v$. However, $v'$ is in the orbit of $w_+$, so $v'$ is not stable, and so $P(v') = 0$. But this implies $\overline P ( \overline v) = \overline{P(v')}= 0$, contradicting our choice of ${\ov v}$. 
	\qed

	\begin{Lemma} \label{Lemma-2}
		Suppose $v \in \ZV_\cO$ is not stable for the action of ${\ZZG}_{\Kal}$ on $\ZV_{\Kal}$. Then the image $\overline v$ of $v$ in $\ZV_{\fF}$ is not stable for the action of $\ZZG_{\fF}$. In particular, if $\ZV_{\fF}$ contains stable vectors, then $\ZV_{\Kal}$ contains stable vectors.
	\end{Lemma}

	\textbf{Proof.}
	Since $v \in \ZV_\cO$ is not stable as an element of $\ZV_{\Kal}$, there exists a nontrivial one-parameter subgroup $\gl: \mathbb{G}_m \to \ZZG_{\Kal}$ such that $I_{\Kal}(\gl, v) = \varnothing$.
	Consider the parabolic subgroup $P_\lam$ of $\ZZG_{\Kal}$ defined by
	\begin{eqnarray*}
		P_\gl({\Kal}) &=& \{ g \in \ZZG_{\Kal}({\Kal}) \mid \lim_{t \to 0} \lam(t)g\lam(t)^{-1} \text{exists in } \ZZG_{\Kal}({\Kal}) \} \\
		&=& \{ g \in \ZZG({\Kal})  \mid gv_i \in \bigoplus_{i \leq j} (\ZV_{\Kal})_j \text{ for all } v_i \in (\ZV_{\Kal})_{i} \}.
	\end{eqnarray*}
	Choose $g \in \ZZG({\Kal})$ such that $g\lam g^{-1}$ has image in $\ZZT$. Using the Iwasawa decomposition, we write $g = qp$ where $q \in \ZZG(\Oal)$ and $p \in P_\gl({\Kal})$, and we set $\lam_{{1}} = p\lam p^{-1}$. Note that $I_{\Kal}(\lam_{{1}}, v) = I_{\Kal}(\lam, {p^{-1}} v) = \varnothing$, because by the definition of $P_\gl$, the weights of ${p^{-1}} v$ with respect to $\lam$ can only be larger than the weights of $v$ with respect to $\lam$.
	Set $\mu := g\lam g^{-1} = q \gl_{{1}} q^{-1}$.
	By the identification of cocharacters $\Hom_{\Kal}(\mathbb{G}_m, \ZZT_{\Kal})\simeq \Hom_\Oal(\mathbb{G}_m,\ZZT) \simeq \Hom_{\fF}(\mathbb{G}_m, \ZZT_{\fF})$, we obtain a nontrivial one-parameter subgroup $\overline{\mu}: \mathbb{G}_m \to \ZZG_{\fF}$ corresponding to $\mu \in \Hom_{\Kal}(\mathbb{G}_m, \ZZT_{\Kal})$.

	Next note that for all $w \in \ZV_{\Oal}$, $I_{\fF}(\overline\mu, \overline w) \subset I_{\Kal}(\mu, w)$;
	in particular $I_{\fF}(\overline \mu, \overline{qv}) \subset I_{{\Kal}}(\mu, qv) = I_{\Kal}(\gl_{{1}}, v) = \varnothing$. Hence $\overline{qv}$ is not stable. Since $\overline{qv} = \overline{q}\cdot\overline v$ is in the orbit of $\overline v$, we see that $\overline v$ is not stable. \qed
	
	\begin{Rem} \label{Rem-MP-unstable}
		The above proof was inspired by the proof of \cite[Proposition~4.3]{MP1}. Moy and Prasad show that if $v\in \ZV_\cO$ is an {\textit{unstable}}, i.e. not semistable, vector of $\ZV_{\Kal}$ under the action of $\ZZG_{\Kal}$, then the image $\overline v$ of $v$ in $\ZV_{\fF}$ is unstable under the action of $\ZZG_{\fF}$.
	\end{Rem}
	
	Combining Lemma \ref{Lemma-1} and Lemma \ref{Lemma-2}, we obtain
	
	\begin{Cor} \label{Cor-stable}
		The representation $\ZV_{\Kal}$ contains stable vectors for the action of $\ZZG_{\Kal}$ if and only if the representation $\ZV_{\fF}$ contains stable vectors for the action of $\ZZG_{\fF}$.
	\end{Cor}
	
	The same statement holds for semistability as follows.
	
	\begin{Lemma} \label{Lemma-semistable}
		The representation $\ZV_{\Kal}$ contains semistable vectors for the action of $\ZZG_{\Kal}$ if and only if the representation $\ZV_{\fF}$ contains semistable vectors for the action of $\ZZG_{\fF}$.
	\end{Lemma}
	
	\textbf{Proof.}
	First suppose that $\ZV_{\Kal}$ contains semistable vectors under $\ZZG_{\Kal}$. As mentioned at the beginning of this section,
	this means that there exists a non-constant (hence nonzero)
$\ZZG_{\Kal}$-invariant homogeneous polynomial on $\ZV_{\Kal}$.
	As in the proof of Lemma \ref{Lemma-1}, the choice of a basis of $\ZV$ yields an identification of $f$ with an element of ${\Kal} [x_1, ..., x_n]$, and we may assume without loss of generality that the minimal valuation of the coefficients of $f$ is zero. 
	Thus we can project $f$ to a non-constant homogeneous $\ZZG({\fF})$-invariant element $\ov f \in \Sym  \check \ZV_{{\fF}}$, 
	and there exists ${\ov v} \in \ZV_{{\fF}}$ such that $\ov f({\ov v}) \neq 0$, i.e. ${\ov v}$ is semistable. 
	
	The converse {follows from} \cite[Proposition~4.3]{MP1} as mentioned in Remark \ref{Rem-MP-unstable}.
	\qed

	\section{Representations coming from Moy--Prasad filtrations} \label{section-Moy--Prasad}
	\subsection{Moy--Prasad filtrations} \label{subsection-filtrations}
	Now we let $G$ be an absolutely simple algebraic group which is defined and split over $K$. Let $x$ be a point in the Bruhat--Tits building $\sB(G,K)$ of $G$ over $K$. Then $x$ is contained in some apartment $\sA$ corresponding to a maximal $K$-split torus $T$ of $G$.  Let $\left(X= \Hom_K(T,\bG_m), \Phi, \check X, \check \Phi\right)$ be the root datum of $G$, and let $\Delta = \{ \al_1, \al_2, ..., \al_\ell \} $ denote the set of simple roots of $T$ corresponding to some Borel subgroup containing $T$. For each $\al \in \Phi$, let $U_\al$ be the unipotent subgroup of $G$ normalized by $T$ and corresponding to $\al$.
	We fix a Chevalley system $\{u_\alpha: \bG_a \to U_\al\}_{\alpha \in \Phi}$, which defines a Chevalley basis $\{ e_\al, h_i \mid \al \in \Phi, 1 \leq i \leq \ell \}$ of $\fg = \Lie(G)$ by $e_\al := \Lie(u_\al)(1)$ and $h_i = [e_{\al_i}, e_{-\al_i}]$. This choice yields a hyperspecial point $x_0 \in \sA$, i.e. the unique point in $\sA$ fixed by $u_\al(-1)u_{-\al}(1)u_\al(-1)$ for all $\al \in \Phi$.
	Taking this point as the origin, we can identify the ($\bR \otimes \check X$)-torsor $\sA$ with $\bR \otimes \check X$. Elements of $X$ are then regarded as affine functions on $\sA$ vanishing at $x_0$. The affine roots are precisely the affine functions $\psi: \sA \to \mathbb{R}$ of the form $\psi = \al + n$, where $\al \in \Phi$ and $n \in \bZ$. We denote the set of affine roots by $\Root$, and for $\psi \in \Root$ denote by $\dot \psi$ the gradient of $\psi$. {We call $x$ \textit{rational} if for every $\psi \in \Psi$, we have $\psi(x) \in \bQ$. For rational $x$, the \textit{order} of $x$ is the smallest positive integer $m$ such that $\psi(x) \in \frac{1}{m}\bZ$ for all $\psi \in \Psi$.}
		 
		 Given $\psi \in \Root$, we define the compact subgroup $U_\psi$ of the root group $U_{\dot \psi} \subset G(K)$ by
	\begin{equation*} U_\psi = \{ u_{\dot \psi}(b) \mid b \in K, \val(b) \geq \psi(x_0) \}.
	\end{equation*} 
	We let $T(K)_0$ denote the maximal bounded subgroup of $T(K)$. Note that
	\begin{equation*} T(K)_0= \{t \in T(K) \mid \val(\chi(t)) = 0 \text{ for all } \chi \in X \}.
	\end{equation*} 
	For $r > 0$, we define the subgroups $T(K)_r$ of $T(K)_0$ by 
	\begin{equation*}
		T(K)_r = \{t \in T(K)_0 \mid \val(\chi(t)-1) \geq r \text{ for all } \chi \in X \}. 
	\end{equation*}

	Then, for $r \geq 0$, the Moy--Prasad filtration subgroups of $G$ are given by
	\begin{equation*}
		G(K)_{x, r} = \< U_\psi, T(K)_r \mid \psi(x) \geq r \>.
	\end{equation*}
	We also set
	\begin{equation*}
		G(K)_{x, r+} = \underset{s > r}\bigcup G(K)_{x, s}.
	\end{equation*}
	The quotient  $G(K)_{x, 0} / G(K)_{x, 0+}$ forms the ${\fF}$-points of a reductive group {$\RP_x$} over ${\fF}$. 
	If we let $r(x)$ be the smallest positive value in the set $\{ \psi(x) \mid \psi \in \Root \}$, then the quotient $\RV_x : = G(K)_{x, r(x)}/G(K)_{x, r(x)+}$ is an ${\fF}$-vector space that admits a rational action of $\RP_x$ induced by conjugation in $G(K)$. 
	Similarly, one can define a filtration $\fg_{x,r}$ ($r \in \bR$) of the Lie algebra $\fg$ such that $G(K)_{x,r}/G(K)_{x,r+}\simeq \fg_{x,r}/\fg_{x,r+}$ as $\RP_x$-representations. See \cite{MP1, MP2} for details.
	
	Note that a slight variation of the proof of \cite[Lemma 3.1]{ReederYu} shows that if $\check \RV_x$ contains semistable vectors, then $x$ is a barycenter of a facet (\cite{ReederYu} shows that if $G(K)_{x, \frac{1}{m}}/G(K)_{x, \frac{1}{m}+}$ contains semistable vectors, then $x$ is a barycenter, but their proof {still} holds after substituting ``$r(x)$" for ``$\frac{1}{m}$"). 
If $x$ is a barycenter, then $x$ is a rational and $r(x) = \frac{1}{m}$, where $m$ is the order of $x$.
	Thus for the rest of the paper, we assume that $r(x) = \frac{1}{m}$ for some integer $m \geq 1$, and that $x$ is rational of order $m$.
	
	\subsection{A reductive group associated to $x$}
	In order to study the action of $\RP_x$ on $\PV_x$ using the results of Section \ref{section-lifting-descending},
	we define an auxiliary reductive group as follows.
	Let $\Phi_x$  be the {sub-root system}  $\{ \al \in \Phi \mid \al(x) \in \bZ \}$ of $\Phi$, and define $H_x$ to be the reductive subgroup of $G$ whose $K$-points are given by $\<T(K), U_\al(K) \mid \al \in \Phi_x \>$. Note that $H_x$ has root datum $(X, \Phi_x, \check X, \check \Phi_x)$. 
	
	Define $Q_x$ to be the parahoric subgroup of $H_x(K)$ generated by $T(K)_0$ and $\{ U_\psi \mid \psi \in \Root, \psi(x) = 0 \}$. 
	Let $\ZG_x$ be the associated parahoric $\cO$-group scheme defined by \cite{MP2}, whose generic fiber $\ZG_x \times_\cO K$ is $H_x$ and with $\ZG_x(\cO) = Q_x$. By construction, $Q_x$ is hyperspecial in $H_x(K)$, so $\ZG_x$ is a split reductive group over $\cO$ {\cite[3.4.2]{Tits}}. By comparing  root data, we see that the special fiber $\ZG_x \times_\cO {\fF}$ of $\ZG_x$ is isomorphic to $\RP_x$. More precisely, we have an isomorphism $\iota$, which on ${\fF}$-points is induced by the inclusion $\ZG_x(\cO)=Q_x \hookrightarrow G(K)_{x,0}$: 
	\begin{equation} \label{eqn-iso1}
		\iota({\fF}): \ZG_x({\fF}) \simeq \ZG_x(\cO)/\ker(\ZG_x(\cO) \twoheadrightarrow \ZG_x({\fF})) {\overset{\sim}\longrightarrow} G(K)_{x,0}/G(K)_{x,0+} \simeq \RP_x({\fF}).
	\end{equation}
	
	Next we construct an $\cO$-module with an action of $\ZG_x$ such that the action on the special fiber corresponds to the action of $\RP_x$ on $\PV_x$.
	Given $0 < r < 1$, we define the $\cO$-submodule $L_{x, r}$ of $\fg$ to be the free $\cO$-submodule with basis 
	\begin{equation*}\{ e_{\alpha+n}:= \varpi^ne_\alpha \mid \alpha \in \Phi, {n \in \bZ}, \alpha(x)+n = r \}.
	\end{equation*}
	Then $L_{x, r}$ is a direct summand of the Moy--Prasad filtration lattice $\fg_{x, r}$, and the inclusion of $L_{x,r}$ into $\fg_{x,r}$ induces an isomorphism 
	\begin{equation} \label{eqn-iso2}
		L_{x, r} \otimes_\cO {\fF} \simeq L_{x,r}/\varpi L_{x,r} {\overset{\sim}\longrightarrow} \fg_{x, r} / \fg_{x, r +}
	\end{equation}
	of ${\fF}$-modules.

	\begin{Lemma}
		The action of $Q_x$ on $\fg$ by restriction of the adjoint action of $G(K)$ stabilizes the $\cO$-module $L_{x, r}$.
	\end{Lemma}
	\textbf{Proof.}
	Let us consider a basis element $e_\gamma$, where 
	$\gamma \in \Psi$ and $\gamma(x) = r$.
	Note that an element of $T(K)_0$ acts on $e_\gamma$ as multiplication by an element of $\cO$, so $T(K)_0$ preserves $L_{x, r}$. Let  $\psi \in \Root$ with $\psi(x) = 0$. Since there exist
	 integers $M_{\dot\psi, \dot\gamma, i}$ such that for any $t \in K$ we have
			$$ u_{\dot \psi}(t)\cdot e_{\dot \gamma} = \sum_{i \geq 0 \atop i\dot\psi + \dot\gamma \in \Phi} M_{\dot\psi, \dot\gamma, i}t^ie_{i\dot\psi + \dot\gamma}$$
			{(\cite[p. 64]{Carter})}, it is an easy calculation to check that
	$U_\psi \cdot e_\gamma \subset L_{x, r}$.
\qed

	As a corollary the adjoint action of $G$ on $\fg$ restricts to an action of $H_x=\ZG_x \times_{\cO} K$ on $L_{x,r}\otimes K$ with the property that $\ZG_x(\cO)$ preserves $L_{x,r}$. Thus by \cite[Section~1.7]{BT2} the action extends to a unique action of $\ZG_x$ on $L_{x,r}$.
	
	Moreover, by construction, if $r(x) \neq 1$, then the action of $\ZG_x \times {\fF}$ on $L_{x,r(x)} \otimes {\fF}$ corresponds via the isomorphisms in (\ref{eqn-iso1}) and (\ref{eqn-iso2}) to the action of $\RP_x$ on $\PV_x$ described in Section \ref{subsection-filtrations}. If $r(x)=1$, then $x$ is hyperspecial in $\sA$ and the $\RP_x$-representation $\PV_{x}$ is isomorphic to the adjoint representation of $\RP_x$ on $\Lie(\RP_x)$, with $\Lie(\RP_x)\simeq \Lie(\ZG_x)_{{\fF}}$.
	Thus we define
	\begin{equation*}
	\ZV_x := \left\{ \begin{array}{cc} 
		L_{x, r(x)} = L_{x, \frac{1}{m}} & \text{if } r(x) \neq 1 \\
		\Lie(\ZG_x) & \text{ if } r(x)=1 \, .
		\end{array}
		\right.   
	\end{equation*}
	 Using Corollary \ref{Cor-stable} and Lemma \ref{Lemma-semistable} we conclude the following:

	\begin{Lemma} \label{Lemma-stableupanddown}
		Let $H_x$ and $\ZV_x$ be as defined above. The action of $(H_x)_{\Kal}$ on $(\ZV_x)_{\Kal} = \ZV_x \otimes_\cO {\Kal}$ has stable (respectively semistable) vectors if and only if the action of $\RP_x$ on $\RV_x$ has stable (respectively semistable) vectors.
	\end{Lemma}
	\begin{Rem} \label{Rem-stableupanddown}
		By the same reasoning, the representation of $(H_x)_{\Kal}$ on ${{({\check\ZV}_x)_{\Kal}}}$ has (semi)stable vectors if and only if the representation of $\RP_x$ on $\check \RV_x$ does. {Here $\check \ZV_x$ denotes the linear dual of $\ZV_x$ over $\cO$, and $\check \RV_x$ denotes the linear dual of $\RV_x$ over $\ov \bF_p$.} In fact, it is this case we are particularly interested in, since the construction of Reeder and Yu requires a stable vector in $\check \RV_x$.
	\end{Rem}

	\section{Vinberg gradings and stability} \label{section-Vinberg}
	
	In this section we use results in Vinberg theory to classify all points $x$  of $\sB(G,K)$ such that $\check \RV_x$ has stable vectors.
	First we show that the representations of $(H_x)_{\Kal}$ on $({\ZV_x})_{\Kal}$ are those coming from Vinberg theory of graded Lie algebras.  We keep the notation from the previous section. Let $G^{ad}$ be the adjoint group of $G$. We have an isogeny $\varphi: G \ra G^{ad}$, and we denote the maximal torus $\varphi(T)$ by $T^{ad}$. The isogeny $\varphi$ induces an inclusion of cocharacters $\check X \hookrightarrow \check X^{ad}$ that yields an isomorphism $\check X \otimes \bR {\overset{\sim}\longrightarrow} \check X^{ad} \otimes \bR$. Using this isomorphism to identify  $\check X \otimes \bR$ and $\check X^{ad} \otimes \bR$, we can write $x$ as $x_0 - \frac{1}{m}\check\lam$ for some $\check \lam \in \check X^{ad}$, where $x_0$ is the hyperspecial point associated to our Chevalley system and $m$ is the order of $x$, both as defined in Section \ref{subsection-filtrations}. Then $x$ induces a homomorphism 
	\begin{eqnarray*}
		\bG_m \overset{\check\lam}\longrightarrow T^{ad} \overset{\Conjug}\longrightarrow \Aut(G), 
	\end{eqnarray*}
	and the choice of a primitive $m$th root of unity $\zeta$ in ${\Kal}$ yields an automorphism 
	$\theta=\Conjug(\check \lambda (\zeta))$ of $G_{\Kal}$. The corresponding automorphism {$d\theta:={\Lie(\theta)}$} of $\fg_{\Kal}$ yields a grading 
	\begin{eqnarray} \label{eqn-g-grading}
		\fg_{\Kal} = \underset{i \in \bZ/m\bZ}\oplus (\fg_{\Kal})_i,
	\end{eqnarray}
	where $\check\lam(\zeta)\cdot v = \zeta^i v$ for all $v \in (\fg_{\Kal})_i$.
	
	On the previously fixed Chevalley basis $\{e_\al, h_{{i}} \}_{\alpha \in \Phi, {1 \leq i \leq \ell}}$, we can write the action of $d\theta$ explicitly:
	\begin{eqnarray*}
		d\theta \cdot e_\al &=& \zeta^{<\check\lam, \al>}e_\al = \zeta^{-{m}\al(x)}e_\al\\
		d\theta \cdot h_{{i}} &=& h_{{i}}.
	\end{eqnarray*}
	
	We see that $(\fg_{\Kal})_{0} = \Lie(H_x){_{\Kal}}$, and hence $(H_x)_{\Kal}$ is the connected component $({G_{\Kal}^\theta})^0$ of the fixed{-}point subgroup of $G_{\Kal}$ under the action of $\theta$.  Moreover, $(\fg_{\Kal})_{-1} = ({\ZV_x})_{\Kal}$, and thus we have an isomorphism 
	\begin{equation} \label{eqn-iso}
		({\check \ZV_x})_{\Kal} \simeq (\fg_{\Kal})_{1}
	\end{equation}	
	of $(H_x)_{\Kal}$-modules, where the action of $({G_{\Kal}^\theta})^0$ on $(\fg_{\Kal})_i$ arises from the adjoint action of $G$ on $\fg$.
	
	The benefit of relating our set-up to that of graded Lie algebras in characteristic zero is that those $(\fg_{\Kal})_1$ which contain stable vectors have been classified in \cite[Corollary~14]{RLYG}. 
		Recall that an element $w$ in the Weyl group $W$ of $G$ is called \textit{elliptic} if $w$ fixes no nonzero vector in $\check X$, and $w$ is $\bZ$-{\textit{regular}} if the group generated by $w$ acts freely on $\Phi$. 	
				
	\begin{Prop}  \label{Prop-principal}
		There exist stable vectors in ${(\check \ZV_x)_{\Kal}}$ under the action of $(H_x)_{\Kal}$ if and only if $m$ is the order of an elliptic $\bZ$-regular element of the Weyl group $W$ of $G$ and $x$ is conjugate under the affine Weyl group $\Waff$ to $x_0 + \frac{1}{m}\check\rho$.
	\end{Prop}
		
	\textbf{Proof.} By \cite[Corollary~14]{RLYG}, $(\fg_{\Kal})_1$ contains stable vectors if and only if $W$ contains {an elliptic $\bZ$-regular element of order $m$} and $d\theta$ is principal, i.e. $\check\lam(\zeta)$ is conjugate by an element of $G^{ad}(\Kal)$ to $\check\rho(\omega)$ for some primitive $m$th root of unity $\omega$. {(Note that \cite{RLYG} formulates the result only for adjoint groups, but one can easily deduce that it holds in our setting as well.)} In addition, {the existence of an elliptic $\bZ$-regular element of order $m$ in $W$ implies that} $\check \rho(\omega)$ is conjugate to $\check\rho(\zeta)$ (\cite[Proposition~8]{RLYG}). But {it is an easy exercise to check that} $\check\lam(\zeta)$ is conjugate to $\check\rho(\zeta)$ if and only if $x$ is conjugate to $x_0 + \frac{1}{m}\check \rho$ under the extended affine Weyl group.
	
	To see that in this case $x$ is conjugate to $x_0 + \frac{1}{m}\check\rho$ under the (unextended) affine Weyl group $\Waff$, we
	note that every element of the extended affine Weyl group can be written as a product $ws$, where $w \in \Waff$ and $s$ is in the stabilizer of an alcove whose closure contains $x_0 + \frac{1}{m}\check \rho$. Checking the normalized Kac Coordinates of all such $x_0 + \frac{1}{m}\check \rho$ (see \cite{RLYG} and \cite{ReederYu}), one verifies that $s$ fixes $x_0 + \frac{1}{m}\check \rho$, and hence $x_0 + \frac{1}{m}\check \rho$ and $x_0-\frac{1}{m}\check\lam$ are conjugated under $\Waff$. \qed

	\subsection{Stable vectors}\label{subsection-stablevectors}
	
	With all the pieces in place, we now come to the main purpose of our paper: the classification of points $x$ such that $\check \RV_x$ contains stable vectors for the action of $\RP_x$.

	\begin{Thm} \label{main-thm} Let $x \in \sA$. Then the representation $\check \RV_{x}$ contains stable vectors under the action of $\RP_x$ if and only if $x$ is conjugate under the affine Weyl group $\Waff$ to $x_0 + \frac{1}{m}\check\rho$, where $m$ is the order of an elliptic, $\bZ$-regular element in the Weyl group $W$ of $G$.	\end{Thm}
	
	\textbf{Proof.}
	This follows directly from Remark \ref{Rem-stableupanddown} and Proposition \ref{Prop-principal}. \qed
	
	 We would like to remark that Reeder and Yu have already given a proof of Theorem \ref{main-thm} for the case {in which} the characteristic $p$ {of the residue field} is sufficiently large {(}see \cite[Corollary~5.1]{ReederYu}{)}. 
	Because
	 it uses Vinberg theory in characteristic $p$, their proof does not {hold} for some small primes, in particular, primes that divide $m$. 
	
	Theorem \ref{main-thm} allows us to use the construction of \cite{ReederYu} to form supercuspidal representations in a uniform way for all $p$. We now briefly review the construction. For details, see \cite[Section 2.5]{ReederYu}. 
	{ Let $\Gk$ be an absolutely simple split reductive group over {the local field} $k$, and let $ = \bF_q$ be the residue field of $k$.  We set $G=\Gk_K$, and let $x {\in} \sB(G,K)$ {be a rational point} that is fixed under the action of the Galois group $\Gal(K/k)$. By identifying the Bruhat--Tits building $\sB(\Gk,k)$ of $\Gk$ over $k$ with the $\Gal(K/k)$-fixed points of $\sB(G,K)$, we may view $x$ as a point of $\sB(\Gk,k)$. Then $G(K)_{x,r}$ is $\Gal(K/k)$-stable {for all $r$;} we denote its $\Gal(K/k)$-fixed points ${(G(K)_{x,r})}^{\Gal(K/k)}$ by $\Gk(k)_{x,r}$. Moreover, the action of $\RP_x$ on $\PV_{x}$ is defined over $\ff$ with $\RP_x(\ff)=\Gk(k)_{x, 0}/\Gk(k)_{x, r(x)}$ and $\PV_{x}(\ff)=\Gk(k)_{x, r(x)}/\Gk(k)_{x, r(x)+}$.}
	 We call a vector $\lam \in \check \PV_x(\ff)$ \textit{stable} if it is stable as a vector in $\check \PV_x$ under the action of $\RP_x$.  
	
	Given a stable vector  $\lam \in \check {{\PV}}_x{(\ff)}$ and a nontrivial character $\chi: \ff^+ \to \bC^\times$, we consider the composition $\chi \circ \lam: \check {{\PV}}_x{(\ff)} \to \bC^\times$ as a character of ${\Gk}(k)_{x, r(x)}$ that is trivial on $\Gk(k)_{x,r(x)+}$. Then the compactly induced representation
	\begin{equation*}
	\pi_x(\lam) := \ind_{{{\Gk}}(k)_{x, r(x)}}^{{\Gk}(k)} (\chi \circ \lam) 
	\end{equation*}
	is a direct sum of irreducible supercuspidal representations of ${\Gk}(k)$ of depth $r(x)$ (\cite[Proposition 2.4]{ReederYu}). 
	
	\begin{Cor}\label{Cor-main}
	Let ${\Gk}$ be a{n absolutely simple} split reductive group over the local field $k$. Then for each $m$ satisfying the conditions of Theorem \ref{main-thm}, there exists a finite unramified extension $k'$ of $k$ such that one can implement the construction of \cite{ReederYu} to form supercuspidal epipelagic representations of ${\Gk}(k')$.
	\end{Cor}

	\textbf{Proof.} Let $G = {\Gk}_K$ {as above}, let $\sA$ be an apartment of $\sB(G,K)$ corresponding to a $k$-split maximal torus of $G$, i.e. $\sA \subset \sB(G,K)^{\Gal(K/k)}$, 
	and let $x = x_0 + \frac{1}{m}\check \rho$ for some $m$ satisfying the conditions of Theorem \ref{main-thm}. By Theorem \ref{main-thm}, the representation of $\RP_x$ on $\check \RV_x$ contains a stable vector, call it $\lam$. Since $\RP_x$ and $\check \RV_x$ are defined over  $\mathfrak{f} = \bF_{q}$, 
	we have that $\lam \in \check \RV_x(\mathfrak{f}')$ for some finite extension $\mathfrak{f}'$ of $\mathfrak{f}$. Let $k'$ be any finite extension of $k$ whose residue field is a finite extension of $\mathfrak{f}'$.
			Then
	 we can input $\lam$ into the construction described above to form $\pi_x(\lam)$, which {decomposes into} a direct sum of
supercuspidal epipelagic representations of ${\Gk}(k')$. 
	 \qed

	\begin{Rem} \label{Rem-main} { Corollary \ref{Cor-main} provides previously unknown representations for some small primes $p$ depending on the type of $\Gk$.} In particular, the construction of \cite{Yu} requires a tameness assumption which fails for small $p$. \end{Rem}

	\subsection{Semistable vectors}
	In \cite[Theorem 8.3]{ReederYu} Reeder and Yu use results in Vinberg theory to classify those $x$ for which $\check \RV_x$ has semistable vectors in terms of conditions that are independent {of} the residue field characteristic $p$. 
	 {Yet} their proof only holds for $p$ sufficiently large: specifically, they assume $p$ is larger than the Coxeter number of $G$.  
	Let us call $x$ a \textit{semistability point} if the prime-independent conditions in \cite[Theorem 8.3]{ReederYu} are satisfied. The proof {of the theorem} can be applied to the setting of graded Lie algebras in characteristic zero to show that $(\check \ZV_{x})_{{\Kal}} \simeq (\fg_{{\Kal}})_{1}$ contains semistable vectors under $(H_x)_{\Kal} = (G_{\Kal}^\theta)^0$ if and only if $x$ is a semistability point. This allows us to extend their characterization to all primes $p$ in the case when the number $r$ appearing in \cite[Theorem 8.3]{ReederYu} is equal to $r(x)$.

	\begin{Prop} \label{Prop-semistable} Let $x \in \sA$ be a rational point of order $m$. 
	Then the representation $\check \RV_{x}$ contains semistable vectors under the action of $\RP_x$ if and only if $x$ is a semistability point.
	\end{Prop}
	
	\textbf{Proof.} The Proposition follows from  Remark \ref{Rem-stableupanddown} using the isomorphism (\ref{eqn-iso}) discussed above. \qed

\section{Classifying stable vectors: an example for $G_2$} \label{section-G2}

Theorem \ref{main-thm} gives necessary and sufficient conditions for the stable locus of $\RP_x$ in ${\check\RV}_x$ to be nonempty. However, determining the stable locus itself 
is currently an ad hoc process that depends on an explicit realization of the representation of $\RP_x$ on ${\check\RV}_x$ in each case.
{Below} we find the stable locus in the case when $G = G_2$ and $x = x_0 + \frac{1}{2}\check\rho$. Up to conjugation, this choice of $x$ determines a unique representation of $G_2(\bQ_2)$, for which we describe the Langlands parameter associated by the local Langlands correspondence. 
This Langlands parameter has characteristics which distinguish it from those corresponding to epipelagic representations of $G(k)$ when $p$ is large.
		
We start with a variation of the Hilbert--Mumford Criterion. Let ${H}$ be a reductive group over an algebraically closed field $E$, and let $V$ be a rational representation of ${H}$ over $E$. Fix a maximal torus in ${H}$, and let $\check X$ denote its cocharacter group. Let $V_s$ be the set of stable vectors in $V$. Recall that $I_E(\lam, v)$ is the set of negative weights for $v$ with respect to $\lam$ (see Section \ref{section-lifting-descending} for details). 
\begin{Lemma}\label{Lemma-hilmum}
Suppose $Y \subset V$ is a ${H}$-invariant subset with the property that $I_E(\mu, y) \neq \varnothing$ for all nontrivial $\mu \in \check X, y \in Y$. Then $Y \subset V_s$. 
\end{Lemma}
\textbf{Proof.}
Suppose $y \in Y$, and let $\gl: \mathbb{G}_m \to {H}$ be a nontrivial one-parameter subgroup. 
	Then there exists $g \in {H}(E)$ such that $\mu : = g\gl g^{-1} \in \check X$. 
By ${H}$-invariance of $Y$, $gy \in Y$, so $I_E(\lam, y) = I_E(\mu, gy) \neq \varnothing$, and $y \in V_s$ by the usual Hilbert--Mumford Criterion.
\qed

\begin{Rem}\label{Rem-invarpoly}
We will apply Lemma \ref{Lemma-hilmum} in the following setting: Let $f$ be a ${H}$-invariant polynomial on $V$ with the property that $f(v) = 0$ whenever
there exists $\mu \in \check X$ such that $I_E(\mu, v) = \varnothing$. This means that $f$ vanishes on vectors which are not stable for the action of the fixed maximal torus. By the lemma, if $y \in V$ and $f(y) \neq 0$, then $y$ is a stable vector for the action of ${H}$. Thus the problem of finding stable vectors is reduced to finding such a polynomial $f$. 
\end{Rem}

For the rest of the section, we take the group $G$ to be {the split form of} $G_2$ {defined over $K$}. Applying Theorem \ref{main-thm} to this case, we see that $\check \PV_x$ has stable vectors under the action of $\RP_x$ if and only if $x$ is conjugate under the affine Weyl group to $x_0+\frac{1}{2}\check \rho$, $x_0+\frac{1}{3}\check \rho$, or $x_0+\frac{1}{6}\check \rho$. When $x = x_0 + \frac{1}{6}\check\rho$, we have that $\RP_x$ is a torus, and the stable vectors in $\check \RV_x$ are easily classified  \cite[Section 2.6]{ReederYu}). For $x = x_0 + \frac{1}{3}\check\rho$, the stable vectors in $\check \RV_x$ are given in \cite[Section 7.5]{ReederYu} in the case $p \neq 3$, and it is not hard to see that this classification extends to the case when $p = 3$. 
		In this section, we will characterize all stable vectors when $x=x_0+\frac{1}{2}\check \rho$.

Fix $x=x_0+\frac{1}{2}\check \rho$. 
We have that $\RP_x \simeq \SO_4 \simeq (\SL_2 \times \SL_2) / \mmu_2$, where $\mmu_2$ is the diagonally embedded group of second 
roots of unity over $\ov \bF_p$. For any commutative ring $A$, let $P_n(A)$ be the space of homogeneous degree-$n$ polynomials over $A$ in two variables, with natural action of $\SL_2(A)$ by precomposition by the transpose, e.g. 
$g = \begin{pmatrix} a & b\\ c & d \end{pmatrix}$ acts by $g \cdot f(X, Y) = f(aX + cY, bX + dY)$. 
 Let $(P_n \boxtimes P_{{l}})(A)$ denote the $A$-span of the standard basis vectors $\{Z^jW^{n-j} \otimes X^kY^{{l}-k} \}_{0 \leq j \leq n, 0 \leq k \leq {l}}$, with natural action of $\SL_2(A) \times \SL_2(A)$, which induces an action of $((\SL_2 \times \SL_2)/\mmu_2)(A)$. By explicitly calculating the action of generators of $\ZG_x(\cO)$ on basis vectors $e_\gamma$, one can check that as an $\ZG_x$-module, $\check L_{x, \frac{1}{2}}$ is isomorphic to $(P_1 \boxtimes P_3)(\cO)$, and as a $\RP_x$-module, $\check\RV_x \simeq (P_1 \boxtimes P_3)({\fF})$.

\subsection{An invariant of $P_1 {\boxtimes} P_3$}\label{subsect-poly}

To use the strategy outlined in Remark \ref{Rem-invarpoly}, we will construct a homogeneous polynomial on $\check \RV_x$ invariant under the action of $\RP_x$. Using the isomorphism of $\cO$-modules $\check L_{x, \frac{1}{2}} \simeq (P_1 \boxtimes P_3)(\cO)$ mentioned above, we see that it suffices to construct a polynomial on $(P_1 \otimes P_3)(\cO)$ invariant under the action of $\ZG_x(\cO) \simeq  (\SL_2 \times \SL_2)/{\mmu}_2(\cO)$. Working over integers
initially will allow us to find a polynomial which works uniformly for all $p$.

First we consider an arbitrary element
\begin{eqnarray*}
F(X, Y) = (aZ + bW) \otimes X^3 + (cZ + dW) \otimes X^2Y + (eZ + fW) \otimes XY^2 + (gZ + hW) \otimes Y^3
\end{eqnarray*}
of ${(}P_1 {\boxtimes} P_3{)}(\cO)$ with $a, b, c, d, e, f, g, h \in \cO$ as a homogeneous degree-{three} polynomial in the variables $X, Y$ with coefficients in $P_1(\cO)$. Then the cubic discriminant $\disc_{X, Y}$ may be written as $\disc_{X, Y} = AZ^4 + BZ^3W + CZ^2W^2 + DZW^3 + EW^4$, where 
	\begin{eqnarray*} 
	A &=& c^2 e^2 - 4 a e^3 - 4 c^3 g + 18 a c e g - 27 a^2 g^2\\
B &=& 2 c d e^2 - 4 b e^3 + 2 c^2 e f - 12 a e^2 f - 12 c^2 d g + 
 18 b c e g + 18 a d e g + 18 a c f g - 54 a b g^2\\ 
 & & - 4 c^3 h + 
 18 a c e h - 54 a^2 g h\\
C &=& d^2 e^2 + 4 c d e f - 12 b e^2 f + c^2 f^2 - 12 a e f^2 - 
 12 c d^2 g + 18 b d e g + 18 b c f g + 18 a d f g\\ 
 & & - 27 b^2 g^2 - 
 12 c^2 d h + 18 b c e h + 18 a d e h + 18 a c f h - 108 a b g h - 
 27 a^2 h^2\\
D &=& 2 d^2 e f + 2 c d f^2 - 12 b e f^2 - 4 a f^3 - 4 d^3 g + 18 b d f g - 
 12 c d^2 h + 18 b d e h + 18 b c f h \\
 &  & + 18 a d f h - 54 b^2 g h -  54 a b h^2 \\
E &=& d^2 f^2 - 4 b f^3 - 4 d^3 h + 18 b d f h - 27 b^2 h^2 .
\end{eqnarray*}
The polynomial $\disc_{X, Y}$ is invariant under the action of the second factor in $\SL_2{(\cO)} \times \SL_2{(\cO)}$. The first factor acts on $\disc_{X, Y}(F)$ via the usual action of $\SL_2$ on {$P_4$}.

Taking the discriminant of the quartic $\disc_{X, Y}(F)$, we may form the composite discriminant $\disc_{Z, W}(\disc_{X, Y} F)$. 
We set $\Delta (F) : = \frac{1}{2^8}\disc_{Z, W}(\disc_{X, Y} F)$, which is a homogeneous polynomial in the variables $a, b, c, d, e, f, g, h$ that turns out to be a polynomial over $\bZ$. 
We can write a formula for $\Delta$ over $\bZ$ as follows. Define two degree-six homogeneous polynomials on $P_1 \boxtimes P_3$ by
\begin{eqnarray*}
H_6(F) &=& -d^3 e^3 + 3 c d^2 e^2 f - 3 c^2 d e f^2 + c^3 f^3 + 9 b d^2 e^2 g + 
 9 b c d e f g - 27 a d^2 e f g - 27 b^2 e^2 f g\\& & - 18 b c^2 f^2 g + 
 27 a c d f^2 g + 54 a b e f^2 g - 27 a^2 f^3 g - 27 b c d^2 g^2 + 
 27 a d^3 g^2 + 81 b^2 c f g^2\\& & - 81 a b d f g^2 - 27 b c d e^2 h + 
 18 a d^2 e^2 h + 27 b^2 e^3 h + 27 b c^2 e f h - 9 a c d e f h - 
 54 a b e^2 f h\\& & - 9 a c^2 f^2 h + 27 a^2 e f^2 h + 54 b c^2 d g h - 
 54 a c d^2 g h - 81 b^2 c e g h + 81 a b d e g h - 81 a b c f g h\\& & + 
 81 a^2 d f g h - 27 b c^3 h^2 + 27 a c^2 d h^2 + 81 a b c e h^2 - 
 81 a^2 d e h^2
 \end{eqnarray*}
and
\begin{eqnarray*}
G_6(F) &=& b c d e f g - a d^2 e f g - b^2 e^2 f g - b c^2 f^2 g + a c d f^2 g + 
 2 a b e f^2 g - a^2 f^3 g - b c d^2 g^2 + a d^3 g^2\\& & + b^2 d e g^2
  +
 2 b^2 c f g^2 - 3 a b d f g^2 - b^3 g^3 - b c d e^2 h + a d^2 e^2 h +
  b^2 e^3 h + b c^2 e f h - a c d e f h\\& & - 2 a b e^2 f h + 
 a^2 e f^2 h + 2 b c^2 d g h - 2 a c d^2 g h - 3 b^2 c e g h + 
 a b d e g h - a b c f g h + 3 a^2 d f g h\\& & + 3 a b^2 g^2 h - 
 b c^3 h^2 + a c^2 d h^2 + 3 a b c e h^2 - 2 a^2 d e h^2 - 
 a^2 c f h^2 - 3 a^2 b g h^2 + a^3 h^3.
 \end{eqnarray*}
Then one can check that $\Delta(F) = H_6(F)^3G_6(F)$. (In fact, $H_6$ and $G_6$ are both invariant polynomials for the action of $\SL_2 \times \SL_2$ on $P_1 \boxtimes P_3$, though we won't need this fact here.)
Let $\overline\Delta$ be {the} reduction {of $\Delta$} mod {$p$}. It is easy to check that $\overline\Delta$ is a nonzero polynomial over $\fF$ on $\check \RV_x$ that is invariant under the action of $\RP_x$.

\subsection{Characterization of Stable Vectors}
\begin{Prop}\label{prop-G2stable}
A vector $F \in {(}P_1 {\boxtimes} P_3{)}({\fF})$ is stable for the action of $\SL_2({\fF}) \times \SL_2({\fF})$ if and only if $\overline \Delta({F}) \neq 0$.
\end{Prop}
\textbf{Proof.}
{Similar to} above, we write an arbitrary vector $F \in (P_1 {\boxtimes} P_3)({\fF})$ as $F = (aZ + bW) \otimes X^3 + {(cZ + dW) \otimes X^2Y + (eZ + fW) \otimes XY^2} + (gZ + hW) \otimes Y^3$, but here we are taking the coefficients $a, b, c, d, e, f, g, h$ 
in ${\fF}$.
Let $\lam (t) = 
\begin{pmatrix}
t^{{s}}& 0\\
0 & t^{-{s}}
\end{pmatrix}
\times
\begin{pmatrix}
t^{{r}} & 0\\
0 & t^{-{r}}
\end{pmatrix}
$ 
be an arbitrary nontrivial cocharacter of the diagonal {maximal} torus. 
 Note that $\lam(t)$ acts on $F$ with weights $3{r} + {s}, 3{r} - {s}, {r} + {s}, {r} - {s}, -{r} + {s}, -{r} - {s}, -3{r} + {s}, -3{r} - {s}$ corresponding to $a, b, c, d, e, f, g, h$.

First we show that if $\overline \Delta(F) \neq 0$, then $F$ is stable.  Assume that $I_{\fF}(\lam, F) = \varnothing$. 
First suppose ${r} + {s} = 0$ and $3{r} - {s} > 0$. Then the weights corresponding to $a, b, c, d, e, f, g, h$ are $2{r}, 4{r}, 0, 2{r}, -2{r}, 0, -4{r}, -2{r}$ respectively, and $I_{\fF}(\lam, F) = \varnothing$ implies $e = g = h = 0$. Thus $G_6(F) = 0$ and so $\overline\Delta(F) = 0$. Similarly, by considering all remaining possible cases given by taking ${r} + {s}$ and $3{r} - {s}$ to be positive, negative, or zero, one obtains that $\overline\Delta(F) = 0$, a contradiction.
Thus $I_{\fF}(\lam, F) \neq \varnothing$, and
using Lemma \ref{Lemma-hilmum} as outlined in Remark \ref{Rem-invarpoly}, we have that if $\overline\Delta({F}) \neq 0$, then ${F}$ is stable.

Next we claim that if $\overline \Delta(F) = 0$, then $F$ is not stable.

First assume $p \neq 2$. 
In this case, $\overline\Delta(F) = 0$ if and only if the polynomial ${\disc_{X, Y}} F$ has a double root, i.e. there is a line{ar form} $l$ in $Z$ and $W$ such that $l^2 \mid \disc_{X, Y} F$. 
Since the first factor of $\SL_2{({\fF})} \times \SL_2{({\fF})}$ acts transitively on {linear forms in} $Z$ {and} $W$, we can assume $Z^2 \mid \disc_{X, Y} F$, i.e. ${\ov D = \ov E} = 0$ {(where $\overline D$ and $\ov E$ are the reductions mod $p$ of the polynomials $D$ and $E$ in Section \ref{subsect-poly})}.
If we write $F$ as $F = Z \otimes F_1(X, Y) + W \otimes F_2(X,Y)$ for cubics $F_1, F_2$, then ${\ov E} = \disc_{X, Y} F_2$.
Since ${\ov E} = 0, F_2(X, Y)$ has a double root, i.e. there is a line{ar form} $l'$ in $X$ and $Y$ such that ${(}l{')}^2 \mid F_2$. Since the second factor of $\SL_2{({\fF})} \times \SL_2{({\fF})}$ acts transitively on line{ar forms} in $X$ {and} $Y$, we can assume $X^2 \mid F_2$, i.e. $f = h = 0$.
Considering the formula for ${\ov D}$, {we see that $dg = 0$}.
If $d = 0$, taking ${r} = 1, {s} = 3$ in $\lam$ above 
gives a one-parameter subgroup such that the weights for $F$ form a subset of $\{ 6, 4, 2, 0 \}$, {giving} $I_{\fF}(\lam, F) = \varnothing$. 
If $g = 0$, taking ${r} = {s} = 1$ gives a one-parameter subgroup such that $I_{\fF}(\lam, F) = \varnothing$. Thus $F$ is not stable.

Now assume $p = 2$. In this case, we use the fact that every stable orbit contains a vector of the form
\begin{equation*}
aZ \otimes X^3 + W \otimes X^2Y + eZ\otimes XY^2 + Z \otimes Y^3,
\end{equation*}
to show that if 
$\overline \Delta (F) = 0$, one can find a one-parameter subgroup $\lam$ such that $I_{\fF}(\lam, F) = \varnothing$. For detailed calculations see the {second author's} PhD thesis (\cite[Section 2.2]{Romano}).

\qed

\begin{Rem}\label{rem-Q2rep}
Suppose $k = \bQ_2$. Using Proposition \ref{prop-G2stable}, one may check that there exists a unique $\RP_x(\ff)$-orbit of stable vectors in  $\check\PV_x(\ff)$ (see Section \ref{subsection-stablevectors} for notation). It is easy to check that if $\lam \in \check\PV(\ff)$ is in this stable orbit, then the stabilizer $\Stab_{\RP_x(\ff)}(\lam)$ is trivial, and by \cite[Proposition 2.4]{ReederYu} the compactly induced representation $\pi_x(\lam)$ of $G_2(\bQ_2)$ is irreducible. In what follows we will describe the Langlands parameter associated to this representation.
\end{Rem}

\subsection{Langlands Parameters}

Let $\Weil_k$ be the Weil group of $k$, and let $\hat G$ be the $\bC$-points of the connected {reductive} group over $\bC$ whose root datum is dual to that of $\Gk$. {According to t}he conjectural local Langlands correspondence one should be able to attach to each irreducible representation $\pi$ of $\Gk(k)$ obtained through the construction of \cite{ReederYu} a discrete Langlands parameter 
\begin{equation*}
\phi: \Weil_k \times \SL_2(\bC) \to \hat G
\end{equation*}
such that certain properties, in particular the formal degree conjecture of \cite{HII}, of the pair $(\pi, \phi)$ are satisfied (here ``discrete" is as defined in \cite[Section 3.2]{GrossReeder}). 
For $p \nmid \lvert W \rvert$, \cite[Section 7]{ReederYu} provides a template for the construction of such parameters. For $p \nmid 2m$, given certain ``epipelagic" parameters, \cite{Kaletha} constructs $L$-packets consisting of epipelagic representations of depth $\frac{1}{m}$ that satisfy many expected properties. 
However, parameters associated to epipelagic representations are not well understood when $p$ is small, except in a few special cases (see \cite[Section 6.3]{GrossReeder} and \cite[Section 7.5]{ReederYu}). Here we will explicitly describe the Langlands parameter associated to the representation discussed in Remark \ref{rem-Q2rep}. We will see that in this case the formal degree conjecture of \cite{HII} is enough to determine a unique parameter. This parameter is not epipelagic in the sense of \cite[Section 4.1]{Kaletha}, and unlike parameters when $p \nmid \lvert W \rvert$, the image of the wild inertia subgroup of $\Weil_k$ is not contained in a torus of $\hat G$.

For the remainder of the section, we fix $k = \bQ_2$, let $\Gk$ be $G_2$ over $\bQ_2$, and let $\pi_x(\lam)$ be the representation of $G_2(\bQ_2)$ discussed in Remark \ref{rem-Q2rep}. For ease of notation, we let $\Weil = \Weil_{\bQ_2}$.
Suppose 
\begin{equation*}
\phi: \Weil \times \SL_2(\bC) \to G_2(\bC)
\end{equation*}
is a discrete Langlands parameter. By \cite[Lemma 3.1]{GrossReeder}, the image $\phi(\Weil)$  of $\Weil$ is finite, so $\phi$ determines a finite, Galois field extension $L_\phi/\bQ_2$ such that
$L_\phi$ is the fixed field of $\ker \phi\mid_{\Weil}$. We have $D :=  \Gal(L_\phi/\bQ_2) \simeq \phi(\Weil)$, and the Lie algebra $\hat\fg$ of $\hat G = G_2(\bC)$ is a representation of $D$, via the compositon of this isomorphism with the adjoint representation. The action of the lower ramification groups 
\begin{equation*}
D \geq D_0 \geq D_1 \geq ... \geq D_c \gneq D_{c+1}= 1,
\end{equation*}
on $\hat\fg$ determines the Swan conductor
\begin{equation*}
b(\phi) = \sum_{j = 1}^c \frac{\dim \hat \fg - \dim {\hat\fg}^{D_j}}{[D_0: D_j]}
\end{equation*}
of this representation, where ${\hat\fg}^{D_j}$ is the subspace of $\hat\fg$ fixed by $D_j$. Recall that the Swan conductor of any representation of $D$ is an integer; in particular $b(\phi) \in \bZ$.

\begin{Prop}
Up to $G_2(\bC)$-conjugacy, there is a unique parameter $\phi: \Weil_{\bQ_2} \times \SL_2(\bC) \to G_2(\bC)$ such that the pair $(\pi_x(\lam), \phi)$ satisfies the formal degree conjecture of \cite{HII}. The filtration of lower ramification groups of $D = \Gal(L_\phi/\bQ_2)$ is given by
\begin{equation}\label{eqn-ramgrps}
D \gneq D_0 \gneq D_1 = D_2 = D_3 \gneq 1,
\end{equation}
where under the isomorphism $D \simeq \phi(\Weil) \subset G_2(\bC)$, $D_1 = D_2 = D_3$ is the unique (up to conjugacy) subgroup of $G_2(\bC)$ isomorphic to $(\bZ/2\bZ)^3$; $D_0$ is the unique subgroup of order $56$ normalizing $D_1$; and $D$ is the unique subgroup of order $168$ normalizing $D_0$. 
\end{Prop}

\textbf{Proof.}
Suppose the pair $(\pi_x(\lam), \phi)$ satisfies the formal degree conjecture of \cite{HII}. Then,
	from \cite[Section 7.1]{ReederYu} we deduce that the following properties hold:
	\begin{enumerate}
		\item[a)] $\phi$ is trivial on $\SL_2(\bC)$
		\item[b)] 
		$b(\phi) = \dim \RP_x = 6$
		\item[c)]
		The centralizer $C_{G_2(\bC)}(\phi(\Weil))$ of the image of $\phi$ is trivial.
	\end{enumerate}

We identify the Galois group $D := \Gal(L_\phi/\bQ_2)$ with $\phi(\Weil) \subset G_2(\bC)$ in order to determine its lower ramification groups.
First we describe $D_c$. Recall that, up to conjugacy, $G_2(\bC)$ contains exactly 3 elementary abelian 2-groups, one of rank $r$ for each $r \leq 3$ (\cite[Theorem 6.1]{Griess}). The reader may refer to \cite{Griess} to determine the normalizers and centralizers of these subgroups.
If $D_c \simeq \bZ/2\bZ$, then $D \subset N_{G_2(\bC)}(D_c) = C_{G_2(\bC)}(D_c)$, hence $D_c \subset C_{G_2(\bC)}(\phi(\Weil))$, 
contradicting (c). 

Now suppose $D_c$ has rank 2. Then $D_c$ is contained in a maximal torus $S$ of $G_2(\bC)$. Let $\fs$ be the Lie algebra of $S$, and let  
$
\delta = \sum_{i = 1}^c \frac{1}{[D_1:D_j]}.
$
Then 
\begin{equation*}
6 = b(\phi) =\frac{12}{[D_0:D_1]}\delta+ \sum_{j = 1}^c \frac{2 - \fs^{D_j}}{[D_0:D_j]},
\end{equation*}
so $\frac{\delta}{[D_0:D_1]} \leq \frac{1}{2}$. 
Since $D \subset N_{G_2(\bC)}(D_c) = N_{G_2(\bC)}(S)$ stabilizes $\fs$, we have
$b_0 := \sum_{j = 1}^c \frac{2 - \fs^{D_j}}{[D_0:D_j]} \in \bZ$. 
And $b_0 < \frac{2\delta}{[D_0:D_1]}$, so $b_0 = 0$ and $D_1 \subset S$.
But $D_1 \subset S$ implies $\delta \in \bZ$, 
contradicting $b(\phi)=6$ since $2 \nmid [D_0:D_1]$.

Thus $D_c$ must be the unique (up to $G_2(\bC)$-conjugacy) subgroup of $G_2(\bC)$ isomorphic to $(\bZ/2\bZ)^3$.
Since $D_1 \subset C_{G_2{\bC}}(D_c) = D_c$, we must have $D_1 = D_2 = ... = D_c$.
Using the facts 
that $N_{G_2(\bC)}(D_c)/D_c \simeq \GL_3(\bF_2)$ and that $D_0/D_1$ is cyclic of order prime to 2, 
we deduce that $[D_0: D_1] = 7$ and $[D:D_0] = 3$, which by (b) implies $c = 3$.
	Thus the lower ramification groups for $D$ are as given in (\ref{eqn-ramgrps}).

Next we will show that there is a unique finite, Galois field extension $L/\bQ_2$ whose lower ramification groups are given by (\ref{eqn-ramgrps}). Indeed, there is a unique Galois extension $E/\bQ_2$ with residue degree 3 and ramification degree 7. Fix a uniformizer $\pi$ of $E$, let $\cO_E$ be the ring of integers of $E$, and for $j \geq 1$ let $U_j = 1 + \pi^j\cO_E$. We write $\mu_7$ for the group of 7th roots of unity in $E$. Suppose $L/\bQ_2$ is a Galois field extension with ramification groups as in (5). We will show that the norm group $N_{L/E}(L^\times)$ is uniquely determined.

We certainly have $\< \pi \> \times \mu_7 \times U_1^2U_4 \subset N_{L/E}(L^\times)$ (see \cite[Theorem 7.12]{Iwasawa}). Let $M/E$ be the abelian extension of local fields such that $N_{M/E}(M^\times) = \< \pi \> \times \mu_7 \times U_1^2U_4$. Note that $\Gamma := \Gal(E/\bQ_2)$ preserves $N_{M/E}(M^\times)$, so $M/\bQ_2$ is Galois. 
As a $\Gamma$-module, $E^\times/N_{M/E}(M^\times)$ decomposes as a sum of two irreducible submodules, each of order 8 (see \cite[Section 5.3]{Romano}), so there are exactly two nontrivial subextensions $E \subset M' \subset M$ which are Galois over $\bQ_2$.
For one of these extensions $M'$, the lower ramification groups $\Gal(M'/\bQ_2)_j = 1$ for $j \geq 2$. This extension was constructed in \cite[Proposition 6.4]{GrossReeder} and corresponds under the local Langlands correspondence to the unique simple supercuspidal representation of $G_2(\bQ_2)$. The other has lower ramification groups as given in (\ref{eqn-ramgrps}), and so must be $L$.

Since every automorphism of $D$ can be realized as conjugation by an element of $N_{G_2(\bC)}(D)$, we see that properties (a), (b), and (c) determine a unique Langlands parameter $\phi$ up to $G_2(\bC)$-conjugacy.
\qed

\bibliography{Stablevectorbib}

\begin{bibdiv}
\begin{biblist}

\bib{Adler}{article}{
      author={Adler, Jeffrey~D.},
       title={Refined anisotropic {$K$}-types and supercuspidal
  representations},
        date={1998},
        ISSN={0030-8730},
     journal={Pacific J. Math.},
      volume={185},
      number={1},
       pages={1\ndash 32},
}

\bib{BT2}{article}{
      author={Bruhat, F.},
      author={Tits, J.},
       title={Groupes r\'eductifs sur un corps local. {II}. {S}ch\'emas en
  groupes. {E}xistence d'une donn\'ee radicielle valu\'ee},
        date={1984},
        ISSN={0073-8301},
     journal={Inst. Hautes \'Etudes Sci. Publ. Math.},
      number={60},
       pages={197\ndash 376},
}

\bib{Carter}{book}{
      author={Carter, Roger~W.},
       title={Simple groups of {L}ie type},
      series={Wiley Classics Library},
   publisher={John Wiley \& Sons, Inc., New York},
        date={1989},
        ISBN={0-471-50683-4},
        note={Reprint of the 1972 original, A Wiley-Interscience Publication},
      review={\MR{1013112 (90g:20001)}},
}

\bib{GrossReeder}{article}{
      author={Gross, Benedict},
      author={Reeder, Mark},
       title={Arithmetic invariants of discrete {L}anglands parameters},
        date={2010},
        ISSN={0020-9910},
     journal={Duke Math J.},
      volume={154},
      number={3},
       pages={431\ndash 508},
}

\bib{Griess}{article}{
      author={Griess, Robert~L., Jr.},
       title={Elementary abelian {$p$}-subgroups of algebraic groups},
        date={1991},
        ISSN={0046-5755},
     journal={Geom. Dedicata},
      volume={39},
      number={3},
       pages={253\ndash 305},
         url={http://dx.doi.org/10.1007/BF00150757},
      review={\MR{1123145}},
}

\bib{HII}{article}{
      author={Hiraga, Kaoru},
      author={Ichino, Atsushi},
      author={Ikeda, Tamotsu},
       title={Formal degrees and adjoint {$\gamma$}-factors},
        date={2008},
        ISSN={0894-0347},
     journal={J. Amer. Math. Soc.},
      volume={21},
      number={1},
       pages={283\ndash 304},
         url={http://dx.doi.org/10.1090/S0894-0347-07-00567-X},
      review={\MR{2350057}},
}

\bib{Iwasawa}{book}{
      author={Iwasawa, Kenkichi},
       title={Local class field theory},
      series={Oxford Science Publications},
   publisher={The Clarendon Press, Oxford University Press, New York},
        date={1986},
        ISBN={0-19-504030-9},
        note={Oxford Mathematical Monographs},
      review={\MR{863740 (88b:11080)}},
}

\bib{Kaletha}{article}{
      author={Kaletha, Tasho},
       title={Epipelagic {$L$}-packets and rectifying characters},
        date={2015},
        ISSN={0020-9910},
     journal={Invent. Math.},
      volume={202},
      number={1},
       pages={1\ndash 89},
         url={http://dx.doi.org/10.1007/s00222-014-0566-4},
      review={\MR{3402796}},
}

\bib{Kim}{article}{
      author={Kim, Ju-Lee},
       title={Supercuspidal representations: an exhaustion theorem},
        date={2007},
        ISSN={0894-0347},
     journal={J. Amer. Math. Soc.},
      volume={20},
      number={2},
       pages={273\ndash 320 (electronic)},
         url={http://dx.doi.org/10.1090/S0894-0347-06-00544-3},
      review={\MR{2276772 (2008c:22014)}},
}

\bib{MP1}{article}{
      author={Moy, Allen},
      author={Prasad, Gopal},
       title={Unrefined minimal {$K$}-types for {$p$}-adic groups},
        date={1994},
        ISSN={0020-9910},
     journal={Invent. Math.},
      volume={116},
      number={1-3},
       pages={393\ndash 408},
}

\bib{MP2}{article}{
      author={Moy, Allen},
      author={Prasad, Gopal},
       title={Jacquet functors and unrefined minimal {$K$}-types},
        date={1996},
        ISSN={0010-2571},
     journal={Comment. Math. Helv.},
      volume={71},
      number={1},
       pages={98\ndash 121},
}

\bib{Mumfordbook}{book}{
      author={Mumford, David},
       title={Geometric invariant theory},
      series={Ergebnisse der Mathematik und ihrer Grenzgebiete, Neue Folge,
  Band 34},
   publisher={Springer-Verlag, Berlin-New York},
        date={1965},
      review={\MR{0214602}},
}

\bib{Mumford}{article}{
      author={Mumford, David},
       title={Stability of projective varieties},
        date={1977},
        ISSN={0013-8584},
     journal={Enseignement Math. (2)},
      volume={23},
      number={1-2},
       pages={39\ndash 110},
}

\bib{RLYG}{article}{
      author={Reeder, Mark},
      author={Levy, Paul},
      author={Yu, Jiu-Kang},
      author={Gross, Benedict~H.},
       title={Gradings of positive rank on simple {L}ie algebras},
        date={2012},
        ISSN={1083-4362},
     journal={Transform. Groups},
      volume={17},
      number={4},
       pages={1123\ndash 1190},
         url={http://dx.doi.org/10.1007/s00031-012-9196-3},
      review={\MR{3000483}},
}

\bib{Romano}{book}{
      author={Romano, Beth},
       title={On the local {L}anglands correspondence: {N}ew examples from the
  epipelagic zone},
   publisher={ProQuest LLC, Ann Arbor, MI},
        date={2016},
        ISBN={978-1339-68032-3},
  url={http://gateway.proquest.com/openurl?url_ver=Z39.88-2004&rft_val_fmt=info:ofi/fmt:kev:mtx:dissertation&res_dat=xri:pqm&rft_dat=xri:pqdiss:10104498},
        note={Thesis (Ph.D.)--Boston College},
      review={\MR{3517860}},
}

\bib{ReederYu}{article}{
      author={Reeder, Mark},
      author={Yu, Jiu-Kang},
       title={Epipelagic representations and invariant theory},
        date={2014},
        ISSN={0894-0347},
     journal={J. Amer. Math. Soc.},
      volume={27},
      number={2},
       pages={437\ndash 477},
         url={http://dx.doi.org/10.1090/S0894-0347-2013-00780-8},
      review={\MR{3164986}},
}

\bib{Tits}{article}{
      author={Tits, J.},
       title={Reductive groups over local fields},
        date={1979},
     journal={Proceedings of Symposia in Pure Mathematics},
      volume={33},
      number={1},
       pages={29\ndash 69},
}

\bib{Yu}{article}{
      author={Yu, Jiu-Kang},
       title={Construction of tame supercuspidal representations},
        date={2001},
        ISSN={0894-0347},
     journal={J. Amer. Math. Soc.},
      volume={14},
      number={3},
       pages={579\ndash 622 (electronic)},
}

\end{biblist}
\end{bibdiv}

\end{document}